\documentclass[a4paper,reqno,11pt]{amsart}
\usepackage{amsmath, amssymb, amsthm, mathrsfs, enumerate, tikz, BOONDOX-cal, units, enumitem, hyperref}
\usetikzlibrary{positioning, calc, matrix, arrows.meta}
\hypersetup{hidelinks}
\usepackage[margin=0.9in]{geometry}

\def\C*{{\sl C*}-algebra} %C*-algebra
\def\Cs*{{\sl C*}-subalgebra} %C*-subalgebra

\newcommand{\NN}{{\mathbb{N}}}
\newcommand{\setst}[2]{\ensuremath{\left\{{#1}\,|\,{#2} \right\}}} %set with "such that" condition
\newcommand{\Mor}[2]{\ensuremath\textup{Mor}{({#1},{#2})}}
\newcommand{\Morunder}[3]{\ensuremath\textup{Mor}_{#3}{({#1},{#2})}}
\newcommand{\idntyof}[1]{\ensuremath{\textup{id}}_{#1}} % Identity map on...
\newcommand{\Ex}{{\mathcal{E\!x}}}   % category of kernel-cokernel pairs in exact structure denoted by same symbol

\newcommand{\Banone}{{\mathcal{B\!a\!n}^{\lower2pt\hbox{\tiny$1$}}}}
\newcommand{\Baninf}{{\mathcal{B\!a\!n}^{\lower1pt\hbox{\tiny$\infty$}}}}
\newcommand{\BanAinf}{{\mathcal{B\!a\!n}^{\lower1pt\hbox{\tiny$\infty$}}_A}}
\newcommand{\Opone}{{\mathcal{O\!p}^{\lower1.5pt\hbox{\tiny$1$}}}}   % category of operators spaces with complete contractions
\newcommand{\Opinf}{{\mathcal{O\!p}^{\lower.5pt\hbox{\tiny$\infty$}}}}  % category of operator spaces with completely bounded maps

\newcommand{\xrighttail}[1]{{\begin{tikzpicture}[auto] \node(E) {}; \node[right=6mm of E](F) { };  \draw[>->] (E.center) to node {\footnotesize{$#1$}} (F.west);{F} \end{tikzpicture}}\!\!}
\newcommand{\xrighttailvar}[2]{{\begin{tikzpicture}[auto] \node(E) {}; \node[right=#1mm of E](F) { };  \draw[>->] (E.center) to node {\footnotesize{$#2$}} (F.west);{F} \end{tikzpicture}}\!\!}
\newcommand{\xrighttwohead}[1]{{\begin{tikzpicture}[auto] \node(E) {}; \node[right=6mm of E](F) { };  \draw[->>] (E.center) to node {\footnotesize{$#1$}} (F.west);{F} \end{tikzpicture}}\!\!}
\newcommand{\xrighttwoheadvar}[2]{{\begin{tikzpicture}[auto] \node(E) {}; \node[right=#1mm of E](F) { };  \draw[->>] (E.center) to node {\footnotesize{$#2$}} (F.west);{F} \end{tikzpicture}}\!\!}

\newcommand{\kercok}[5]{{#1}\xrighttail{#4}{#2}\xrighttwohead{#5}{#3}} % kernel-cokernel pair of admissible monic and admis epic. Left to right objects then kernel then cokernel
\newcommand{\kercokvar}[6]{{#1}{\xrighttailvar{#6}{#4}}{#2}{\xrighttwoheadvar{#6}{#5}}{#3}} % kernel-cokernel pair of admissible monic and admis epic. Left to right objects then kernel then cokernel

\newcommand{\Midim}[1]{\ensuremath\textup{Inj}_{\mathcal{M}}\textup{-dim}\,{(#1)}} %I-dim
\theoremstyle{definition} %upright "theorems"
\newtheorem{definition}{De{f}inition}[section]
\newtheorem{remark}[definition]{Remark}

\theoremstyle{plain} %all italics
\newtheorem{theorem}[definition]{Theorem}
\newtheorem{prop}[definition]{Proposition}
\newtheorem{cor}[definition]{Corollary}

\let\oldmarginpar\marginpar
\renewcommand\marginpar[1]{\-\oldmarginpar[\raggedleft\bf\scriptsize #1]%
	{\raggedright\bf\scriptsize #1}}

\numberwithin{equation}{section}

\title{Schanuel's Lemma for Exact Categories}
\author{Martin Mathieu}
\address{Mathematical Sciences Research Centre\\
Queen's University Belfast\\
Belfast BT7 1NN\\
Northern Ireland}
\email{m.m@qub.ac.uk}
\author{Michael Rosbotham}
\address{Mathematical Sciences Research Centre\\
Queen's University Belfast\\
Belfast BT7 1NN\\
Northern Ireland}
\email{mrosbotham01@qub.ac.uk}

   % AMS 2020 subject classifications
   % helpful keywords

\begin{document}\bibliographystyle{plain}

\begin{abstract}
We prove an injective version of Schanuel's lemma from homological algebra in the setting of exact categories.
\end{abstract}

	\maketitle
\section{Introduction}\label{sect:intro}

\noindent
Schanuel's lemma is a useful tool in homological algebra and category theory.
It appears to have come about as a response to a question by Kaplansky, see \cite[p.~166]{LamLectures},
and simplifies the definition of the projective (or, injective) homological dimension in module categories,
hence in abelian categories.
The typical categories that arise in functional analysis are not abelian but lately, the use of exact structures
on additive categories of Banach modules and related ones has been suggested and indeed been exploited successfully.

In~\cite{B11}, B\"uhler develops homological algebra for bounded cohomology in the setting of Quillen's exact categories.
In~\cite{AM2}, exact categories of sheaves of operator modules over \textsl{C*}-ringed spaces are studied.
Relative cohomology and cohomological dimension for (not necessarily self-adjoint) operator algebras is the topic
of~\cite{MR2021}, see also~\cite{Rosb2020}. In view of this, it seems beneficial to establish an injective version
of Schanuel's lemma for exact categories and show how it yields the injective dimension theorem.

When we equip an additive category $\mathcal A$ with an exact structure we fix a pair $(\mathcal M,\mathcal P)$
consisting of a class of monomorphisms $\mathcal M$ and a class of epimorphisms $\mathcal P$ such that
each $\mu\in\mathcal M$ and $\pi\in\mathcal P$ form a kernel-cokernel pair which we write as
\[
\begin{tikzpicture}[auto]
	\node(E) {$E$};
	\node[right=11mm of E](F) {$F$};
	\node[right=11mm of F](G) {$G$};
	\draw[>->] (E) to node {\footnotesize{$\mu$}} (F);
	\draw[->>] (F) to node {\footnotesize{$\pi$}} (G);
\end{tikzpicture}
\]
where $E,F$ and $G$ are objects in~$\mathcal A$.
We require that $\mathcal M$ and $\mathcal P$ contain all identity morphisms and are closed under composition,
and term their elements as \textit{admissible monomorphisms\/} and \textit{admissible epimorphisms}, respectively.
Furthermore, the push-out of an admissible monomorphism along an arbitrary morphism exists and yields
an admissible monomorphism, and, likewise, the pull-back of an admissible epimorphism along an arbitrary morphism exists
and yields an admissible epimorphism. If these conditions are fulfilled and $(\mathcal M,\mathcal P)$
is invariant under isomorphisms, $(\mathcal M,\mathcal P)$ is called \textit{an exact structure\/} on~$\mathcal A$
and will typically be denoted by~$\Ex$. The pair $(\mathcal A,\Ex)$ is said to be an
\textit{exact category}.

Unlike in abelian categories not every morphism in an exact category has a canonical factorisation into an epimorphism
followed by a monomorphism. One therefore has to restrict to \textit{admissible morphisms\/} which are those that
arise as $\mu\circ\pi$ for some $\mu\in\mathcal M$ and $\pi\in\mathcal P$. (It is easy to check that, once such
factorisation exists, it is unique up to unique isomorphism.)

The kernel-cokernel pairs replace the usual short exact sequences in abelian categories while long exact sequences
are built from admissible morphisms. A very readable introduction into exact categories is given in~\cite{B08}.

In this note, we provide the details of how Schanuel's lemma works in general exact categories
and establish the Injective Dimension Theorem (Theorem~\ref{Theorem: Injective dimension theorem}).

\section{Preliminaries}\label{Section: Prelim}

\noindent
We include here the necessary terminology and initial results,
for a fixed exact category $(\mathcal A,\Ex)$, where $\Ex=(\mathcal M,\mathcal P)$.

\begin{definition}\label{Definition: M-injective}
An object $I$ in an exact category $(\mathcal{A},\Ex)$
is \textit{$\mathcal{M}$-injective} if, when given $E\xrighttail{\mu}{F}$ and
a morphism $f\in\Morunder{E}{I}{\mathcal{A}}$,  for objects $E,F\in\mathcal{A}$,
there exists a morphism ${g}\in\Morunder{F}{I}{\mathcal{A}}$
making the following diagram commutative
\[
\begin{tikzpicture}[auto]
	\node (E) {$E$};
	\node (F) [right=2cm of E]{$F$};
	\node (I)[below=1.2cm of E]{$I$};
	\draw[>->] (E)--node{\footnotesize{$\mu$}}(F);
	\draw[->] (E)--node[swap]{\footnotesize{$f$}}(I);
	\draw[->, dashed] (F)--node{ }(I);
\end{tikzpicture}
\]
The exact category has \textit{enough $\mathcal{M}$-injectives} if,
for every $E\in\mathcal{A}$, there exist an $\mathcal{M}$-injective object~$I$
and an admissible monomorphism $E\xrighttail{ }{I}$.
\end{definition}

We will also make use of the following characterisations of $\mathcal{M}$-injective objects.

\begin{prop}\label{Prop: M-injective}
Let $E$ be an object in an exact category $(\mathcal{A},\Ex)$.
The following are equivalent.
\begin{enumerate}[label=\upshape(\roman*)]
\item $E$ is $\mathcal{M}$-injective;
\item Every admissible monomorphism $E\xrighttail{ }{F}$, for $F\in\mathcal{A}$,  has a left inverse;
\item There exist an $\mathcal{M}$-injective object $I\in\mathcal{A}$
and a morphism $E\longrightarrow{I}$ with a left inverse (i.e.,~$E$ is a retract of an
$\mathcal{M}$-injective object).
\end{enumerate}
\end{prop}

\noindent
The arguments are standard.

As exact categories are additive,
we can form the product of any two objects (and thus, of any finite number of objects).

\begin{prop}\label{Prop: Product addcats}
Let $E, F, G$ be objects in an additive category~$\mathcal{A}$.
The following are equivalent:
\begin{enumerate}[label=\upshape(\roman*)]
	\item $F$ is a product of $E$ and $G$;
	\item $F$ is a coproduct of $E$ and $G$;
	\item There exist a kernel-cokernel pair in $\mathcal{A}$,
	\begin{equation}\label{Diagram: biproduct}
		\begin{tikzpicture}[auto, baseline=(current  bounding  box.center)]
			\node(E) {$E$};
			\node[right=11mm of E](F) {$F$};
			\node[right=11mm of F](G) {$G$};
			\draw[>->] (E) to node {\footnotesize{$\mu$}} (F);
			\draw[->>] (F) to node {\footnotesize{$\pi$}} (G);
		\end{tikzpicture}
	\end{equation}
	and morphisms $\widetilde{\mu}\in\Morunder{F}{E}{\mathcal{A}}$ and
	$\widetilde{\pi}\in\Morunder{G}{F}{\mathcal{A}}$
	such that
	$\widetilde{\mu}\circ\mu=\idntyof{E}$ and $\pi\circ\widetilde{\pi}=\idntyof{G}$,
	and $\mu\circ\widetilde{\mu}+\widetilde{\pi}\circ\pi=\idntyof{F}$;	
	\item There exist a kernel-cokernel pair in $\mathcal{A}$,
\begin{equation}\label{Diagram: biproduct admon}
\begin{tikzpicture}[auto, baseline=(current  bounding  box.center)]
	\node(E) {$E$};
	\node[right=11mm of E](F) {$F$};
	\node[right=11mm of F](G) {$G$};
	\draw[>->] (E) to node {\footnotesize{$\mu$}} (F);
	\draw[->>] (F) to node {\footnotesize{$\pi$}} (G);
\end{tikzpicture}
\end{equation}
and a morphism $\widetilde{\mu}\in\Morunder{F}{E}{\mathcal{A}}$ such that
$\widetilde{\mu}\circ\mu=\idntyof{E}$, the identity morphism on~$E$.
\end{enumerate}
Moreover, if these equivalent conditions are met,
the kernel-cokernel pair in Diagram~\eqref{Diagram: biproduct admon}
will belong to every exact structure that can be placed on~$\mathcal{A}$.
\begin{proof}
Finite products, coproducts and biproducts coincide in an additive category 
(see, e.g., \cite[Proposition~7.1--Corollary~7.3.]{O2000}),
and condition~(iii) is just the definition of $F$ being a biproduct of $E$ and~$G$.
That condition~(iii) is equivalent to condition~(iv)
can be proven in the exact same way as the `Splitting Lemma' in module theory
(see, e.g., \cite[Proposition~4.3.]{MacLane_Homology}).
The final statement of this proposition
is a direct consequence of the conditions required for monomorphisms and epimorphisms to be admissible;
see \cite[Lemma~2.7.]{B08} for details.
\end{proof}
\end{prop}
Kernel-cokernel pairs satisfying condition (iii) of Proposition~\ref{Prop: Product addcats}
are said to be \textit{split}.
For objects $E$ and $F$ in $\mathcal{A}$ we will denote their (co)product by $E\oplus{F}$.

\begin{prop}\label{Prop: Direct sum kercokpair}
	Suppose $\kercok{E}{F}{G}{\footnotesize{\mu}}{\footnotesize{\pi}}$ is a kernel-cokernel pair in $\Ex$.
	\begin{enumerate}[label=\upshape(\roman*)]
	\item For any $A\in\mathcal{A}$, there is a kernel-cokernel pair in $\Ex$,
	\[
	\kercokvar{E\oplus{A}}{F\oplus{A}}{G}{ }{ }{10}
	\]
	\item If $F\cong E\oplus{G}$, then $F$ is $\mathcal{M}$-injective
	if and only if both $E$ and $G$ are $\mathcal{M}$-injective.
	\end{enumerate}
\begin{proof}
We first prove (i). For $A\in\mathcal{A}$, there exist split kernel-cokernel pairs
\[\begin{tikzpicture}[auto]
	\node(1) {$E$};
	\node[right=9mm of 1] (2) {$E {\oplus} A$};
	\node[right=9mm of 2] (3) {$A$};
	\node[right=10mm of 3] (and) {\textup{and}};
	\node[right=10mm of and](4) {$A$};
	\node[right=9mm of 4] (5) {$F {\oplus} A$};
	\node[right=9mm of 5] (6) {$F$};
	\draw[>->] (1) to node {{\footnotesize{$\iota$}}} (2); \draw[>->] (4) to node {{\footnotesize{$\theta$}}} (5);
	\draw[->>] (2) to node {{\footnotesize{$\rho$}}} (3); \draw[->>] (5) to node {{\footnotesize{$\tau$}}} (6);
\end{tikzpicture}
\]
and $\pi\circ\tau\in\Morunder{F\oplus{A}}{G}{\mathcal{A}}$
is an admissible epimorphism, as a compostion of morphisms in~$\mathcal{P}$.
Define $\varphi\in\Morunder{E\oplus{A}}{F\oplus{A}}{\mathcal{A}}$ by
\[\varphi=\widetilde{\tau}\circ{\mu}\circ\widetilde{\iota} +\theta\circ\rho\]
using the same notation as in Proposition~\ref{Prop: Product addcats}.
Then $(\varphi, \pi\circ\tau)$ is the desired kernel-cokernel pair.

To show this, it is enough to demonstrate that $\varphi$ is a kernel of $\pi\circ\tau$.
First note the composition
\[
(\pi\circ\tau)\circ\varphi =
\idntyof{F}\circ(\pi\circ{\mu})\circ\widetilde{\iota} + \pi\circ(\tau\circ\theta)\circ\rho= 0.
\]
Now suppose there exist $B\in\mathcal{A}$ and a morphism $f\in\Morunder{B}{F\oplus{A}}{\mathcal{A}}$
such that $(\pi\circ\tau)\circ{f}=0$.
As $\mu$ is a kernel for $\pi,$ there exists a unique morphism
$g'\in\Morunder{B}{E}{\mathcal{A}}$ such that
$\mu\circ{g'}=\tau\circ{f}$.
Define $g\in\Morunder{B}{E\oplus{A}}{\mathcal{A}}$ by
\[
g=\iota\circ{g'}+\widetilde{\rho}\circ\widetilde{\theta}\circ{f}.
\]
Then $\varphi\circ{g}=(\widetilde{\tau}\circ\tau)\circ{f} + (\theta\circ\widetilde{\theta})\circ{f}=
\idntyof{F\oplus{A}}\circ{f}=f$.

To finish the proof of (i),
we show that there is no other morphism $h\in\Morunder{B}{E\oplus{A}}{\mathcal{A}}$
such that $\varphi\circ{h}=f$.
Suppose we have such a morphism $h$.
Then,  $\widetilde{\theta}\circ{f} = \widetilde{\theta}\circ\varphi\circ{h} = \rho\circ{h}$,
and $\mu\circ{g'} = \tau\circ f = \tau \circ \varphi\circ h
= \mu \circ \widetilde{\iota}\circ h,$
and therefore $g'=\widetilde{\iota}\circ h$.
Combining these facts gives:
\[
h = \idntyof{E \oplus A}\circ h
 = (\iota \circ \widetilde{\iota} + \widetilde{\rho} \circ \rho)
 =\iota \circ g' + \widetilde{\rho}\circ\widetilde{\theta}\circ{f}
 =g,
\] as required.

For assertion~(ii) suppose $F\cong E\oplus G$. Then there exist morphisms $\widetilde{\mu}\in\Morunder{F}{E}{\mathcal{A}}$ and
$\widetilde{\pi}\in\Morunder{G}{F}{\mathcal{A}}$
such that
$\widetilde{\mu}\circ\mu=\idntyof{E}$ and $\pi\circ\widetilde{\pi}=\idntyof{G}$,
and $\mu\circ\widetilde{\mu}+\widetilde{\pi}\circ\pi=\idntyof{F}$.
In particular, $E$ and $G$ are retracts of~$F$.
By Proposition~\ref{Prop: M-injective}, if $F$ is $\mathcal{M}$-injective so are $E$ and $G$.
Finally, suppose $E$ and $G$ are $\mathcal{M}$-injective and there is an admissible monomorphism
\[\begin{tikzpicture}[auto]
	\node(F) {$F$};
	\node[right=14mm of E](B) {$B$};
	\draw[>->] (F) to node {{\footnotesize{$f$}}} (B);
\end{tikzpicture}
\]
 where $B\in\mathcal{A}$.
 Because $E$ and $G$ are $\mathcal{M}$-injective,
 there exist $g_E\in\Morunder{B}{E}{\mathcal{A}}$
 such that $\widetilde{\mu}=g_E \circ f$
 and $g_G\in\Morunder{B}{G}{\mathcal{A}}$
 such that
 and $\pi= g_F \circ f$.
Let $g= \mu\circ g_E + \widetilde{\pi}\circ g_G$,
then $g$ is a left inverse of $f$, indeed:
 \[ g \circ f=
  \mu\circ (g_E \circ f)  + \widetilde{\pi}\circ (g_G \circ f)
  = \mu\circ\widetilde{\mu}+\widetilde{\pi}\circ\pi=\idntyof{F}.
 \]
 Hence, by Proposition~\ref{Prop: M-injective},
 $F$ is $\mathcal{M}$-injective.
\end{proof}
\end{prop}

\section{Schanuel's Lemma}\label{Section: Schanuel's Lemma}

\noindent
Fix an exact category $(\mathcal{A},\Ex)$. The following is the injective version of Schanuel's lemma for exact categories.

\begin{prop}\label{Prop: Schanuel's Lemma}
Suppose $\kercok{E}{I}{F}{{\footnotesize{\mu}}}{{\footnotesize{\pi}}}$
and  $\kercok{E}{I'}{F'}{{\footnotesize{\mu'}}}{{\footnotesize{\pi'}}}$ are kernel-cokernel pairs in $\Ex$,
and that $I, I'$ are $\mathcal{M}$-injective objects.
Then $I\oplus{F'}\cong I'\oplus{F}$ in $\mathcal{A}$.
\begin{proof}
First, by the axioms of an exact structure, we can form the following push-out,
\begin{equation}\label{Diagram: Square}
	\begin{tikzpicture}[auto, baseline=(current  bounding  box.center)]
		\matrix(m)[column sep=4em, row sep=3em]{
			\node(E) {$E$}; & \node(I) {$I$};
			\\
			\node(I') {$I'$}; & \node(C) {$C$};
			\\};
		\draw[>->] (E) to node {\footnotesize{$\mu$}} (I); \draw[>->] (E) to node [swap]{\footnotesize{$\mu'$}} (I');
		\draw[>->] (I) to node {\footnotesize{$h$}} (C); \draw[>->] (I') to node [swap]{\footnotesize{$h'$}} (C);
	\end{tikzpicture}
\end{equation}
where every morphism is an admissible monomorphism.
Extending this diagram to include the given cokernels, and adding in some zero morphisms, we get the following commutative diagram:
\vspace*{-16mm}
\begin{equation}\label{Diagram: Square with 0s}
\begin{tikzpicture}[auto, baseline=(current  bounding  box.center)]
	\matrix(m)[column sep=4em, row sep=3em]{
		\node(E) {$E$}; & \node(I) {$I$}; & \node(F) {$F$};
		\\
		\node(I') {$I'$}; & \node(C) {$C$}; &{}
		\\
	    \node(F') {$F'$}; & {} &{}\\	};
	\draw[>->] (E) to node {\footnotesize{$\mu$}} (I); \draw[>->] (E) to node [swap]{\footnotesize{$\mu'$}} (I');
	\draw[>->] (I) to node {\footnotesize{$h$}} (C); \draw[>->] (I') to node [swap]{\footnotesize{$h'$}} (C);
	\draw[->>] (I) to node {\footnotesize{$\pi$}} (F); \draw[->>] (I') to node [swap]{\footnotesize{$\pi'$}} (F');
	\draw[->] (I'.south east) to [out=300, in=270, looseness=1.3] node [swap]{\footnotesize{$0$}} (F);
	\draw[->] (I.north) to [out=120, in=180, looseness=2] node [swap]{\footnotesize{$0$}} (F');
\end{tikzpicture}
\end{equation}
By the universal property of push-outs, there are a unique morphism
$p\in\Mor{C}{F}$ such that $ph'=0$ and $ph=\pi$,  and a unique morphism $p'\in\Mor{C}{F'}$
such that $p'h=0$ and $p'h'=\pi'$.
Hence, we have the following commutative diagram:
\begin{equation}\label{Diagram: Square with more squares}
	\begin{tikzpicture}[auto, baseline=(current  bounding  box.center)]
		\matrix(m)[column sep=4em, row sep=3em]{
			\node(E) {$E$}; & \node(I) {$I$}; & \node(F) {$F$};
			\\
			\node(I') {$I'$}; & \node(C) {$C$}; &\node(F2) {$F$};
			\\
			\node(F') {$F'$}; & \node(F'2) {$F'$}; &{}\\	};
		\draw[>->] (E) to node {\footnotesize{$\mu$}} (I); \draw[>->] (E) to node [swap]{\footnotesize{$\mu'$}} (I');
		\draw[>->] (I) to node {\footnotesize{$h$}} (C); \draw[>->] (I') to node [swap]{\footnotesize{$h'$}} (C);
		\draw[->>] (I) to node {\footnotesize{$\pi$}} (F); \draw[->>] (I') to node [swap]{\footnotesize{$\pi'$}} (F');
	\draw[->] (C) to node [swap]{\footnotesize{$p$}} (F2); \draw[->] (C) to node {\footnotesize{$p'$}} (F'2);
	\draw[->] (F) to node {\footnotesize{$\idntyof{F}$}} (F2); \draw[->] (F') to node [swap]{\footnotesize{$\idntyof{F'}$}} (F'2);
	\end{tikzpicture}
\end{equation}
The result will follow if the middle row and middle column are both split kernel-cokernel pairs.
As $h, h'\in\mathcal{M}$ and $I,I'$ are $\mathcal{M}$-injective,
this will be the case if both $(h',p)$ and $(h,p')$ are kernel-cokernel pairs.
We deal with $(h',p)$, the other pair is done in the exact same way.

To show that $(h',p)$ is a kernel-cokernel pair, it is enough to verify that $p$ is a cokernel of~$h'$.
Suppose there exist an object $G\in\mathcal{A}$ and a morphism $q\in\Mor{C}{G}$ such that $qh'=0$.
We are done if we find a unique morphism $\psi\in\Mor{F}{G}$ such that the following diagram is commutative:
\begin{equation}\label{Diagram: cokernel}
\begin{tikzpicture}[auto, baseline=(current  bounding  box.center)]
	\matrix(m)[column sep=4em, row sep=3em]{
\node(I') {$I'$}; & \node(C) {$C$}; &\node(F2) {$F$};
\\
{}; & \node(G) {$G$}; &{}\\	};
\draw[>->] (I') to node [swap]{\footnotesize{$h'$}} (C);
\draw[->] (C) to node [swap]{\footnotesize{$p$}} (F2);
\draw[->] (C) to node [swap]{\footnotesize{$q$}} (G);
\draw[->, dashed] (F2) to node {\footnotesize{$\psi$}} (G);
\draw[->, bend left] (I') to node {\footnotesize{$0$}} (F2);
\draw[->, bend right] (I') to node [swap]{\footnotesize{$0$}} (G);
	\end{tikzpicture}
\end{equation}

We have $(qh)\mu=q(h\mu)=q(h'\mu')=0$ and, because $(\mu,\pi)$ is a kernel-cokernel pair, there exists a unique morphism
$t\in\Mor{F}{G}$ such that $t\pi=qh$.
Therefore, the following diagram is commutative:
\begin{equation}\label{Diagram: Square+G}
	\begin{tikzpicture}[auto, baseline=(current  bounding  box.center)]
		\matrix(m)[column sep=4em, row sep=3em]{
			\node(E) {$E$}; & \node(I) {$I$}; & {}
			\\
			\node(I') {$I'$}; & \node(C) {$C$}; & {}\\
			{}; &  { }; & \node(G) {$G$};
			\\};
		\draw[>->] (E) to node {\footnotesize{$\mu$}} (I); \draw[>->] (E) to node [swap]{\footnotesize{$\mu'$}} (I');
		\draw[>->] (I) to node {\footnotesize{$h$}} (C); \draw[>->] (I') to node [swap]{\footnotesize{$h'$}} (C);
		\draw[->] (C) to node {\footnotesize{$q$}} (G);
		\draw[->, bend right] (I') to node [swap]{\footnotesize{$0$}} (G);
		\draw[->, bend left] (I) to node {\footnotesize{$t\pi$}} (G);
	\end{tikzpicture}
\end{equation}
By the universal property of push-outs,
$q$ is the unique morphism $C\rightarrow{G}$ that makes Diagram~\eqref{Diagram: Square+G} commutative.
However, $(tp)h=t(ph)=t\pi$ and $(tp)h'=t(ph')=0$.
So, $q=tp$ and setting $\psi=t$ makes Diagram~\eqref{Diagram: cokernel} commutative.
Finally, suppose there also exists $t'\in\Mor{F}{G}$ such that $q=t'p$.
Recalling from Diagram~\eqref{Diagram: Square with more squares} that $\pi=ph$, we have
\[ t'\pi= t'(ph)=(t'p)h=(tp)h=t(ph)=t\pi,
\]
and, because $\pi$ is an epimorphism, $t'=t$. Thus, uniqueness has been verified.
\end{proof}
\end{prop}

\begin{cor}\label{Corollary: Schanuel's iso}
Suppose there is a diagram of morphisms in an exact category $(\mathcal A,\Ex)$ of the form
\[
\begin{tikzpicture}[auto]
\matrix[column sep={1.1cm}, row sep={1.3cm,between origins}, matrix of math nodes](m)
{ E{\hphantom{'}} &   I{\hphantom{'}}  & F{\hphantom{'}} \\
  E'&   I'  & F' \\
};
\draw[>->] (m-1-1) to node{ } (m-1-2); \draw[->>] (m-1-2) to node{ } (m-1-3);
\draw[>->] (m-2-1) to node{ } (m-2-2); \draw[->>] (m-2-2) to node{ } (m-2-3);
\draw[->]  (m-1-1) to node[swap]{$\footnotesize{\cong}$} (m-2-1);
\end{tikzpicture}
\]
such that $I$ and $I'$ are $\mathcal{M}$-injective,
the horizontal lines are in $\Ex$ and
the vertical arrow is an isomorphism.
Then $I\oplus{F'}\cong I'\oplus{F}$ in $\mathcal{A}$.
\end{cor}

We extend Schanuel's lemma to injective resolutions in Proposition~\ref{Prop:Schanuel's Lemma resolution} below.
Recall that a morphism is \textit{admissible\/} if it is the composition $\mu\circ\pi$ for some $\mu\in\mathcal M$
and $\pi\in\mathcal P$. Such factorisation is unique up to unique isomorphism (\cite[Lemma~8.4]{B08}).

\begin{definition}\label{Definition: Injective res}
	For an object $E\in\mathcal{A}$,
	an \textit{$\mathcal{M}$-injective resolution} of $E$ is a sequence of admissible morphisms
	of the form:
	\[\begin{tikzpicture}[auto]\matrix[column sep={0.5cm}, row sep={1.3cm,between origins}]{
			\node(E) {$E$}; &{} & \node(0)  {$I^0$};&{} & \node(dots)  {$\cdots$};&{} & \node(1)  {$I^{n-1}$};&{} & \node(m) {$I^n$};&{} & \node(F) {$\cdots$};\\
			\node(not1){ };&\node(E2) {$G^0$};&{}&\node(K1) {$G^1$};&{}&\node(K1m)  {$G^{n-1}$};&{}&\node(Km) {$G^n$}; &{}
			&\node(Km1) {$G^{n+1}$};\\	};	
		
		\draw[>->] (E) to node {} (0); \draw[->] (0) to node {} (dots);
		\draw[->] (dots)to node {} (1); \draw[->] (1) to node {} (m);
		\draw[->] (m) to node {} (F);
		\draw[->>] (E) to node [swap]{{\footnotesize{$\cong$}}} (E2);
		\draw[->>] (0) to node { } (K1);\draw[->>] (dots) to node { } (K1m); \draw[->>] (1) to node { } (Km);
		\draw[>->] (E2) to node { } (0);\draw[>->] (K1) to node { } (dots); \draw[>->] (K1m) to node { } (1); \draw[>->] (Km) to node { } (m);
		\draw[->>]  (m) to node { } (Km1);
		\draw[>->] (Km1) to node { } (F);
	\end{tikzpicture} 	\]
	such that, for each $n\geq0,$ the object $I^n$ is $\mathcal{M}$-injective,
	and
	\[
	\begin{tikzpicture}[auto, baseline=(current  bounding  box.center)]
		\node(n) {$G^n$};
		\node[right=11mm of n](I) {$I^n$};
		\node[right=11mm of I](n+1) {$G^{n+1}$};
		\draw[>->] (n) to node { } (I);
		\draw[->>] (I) to node { } (n+1);
	\end{tikzpicture}
	\]
	forms a kernel-cokernel pair in $\Ex$.

\end{definition}
If $\mathcal{A}$ has enough $\mathcal{M}$-injectives,
we can build an injective resolution for every object in~$\mathcal{A}$.

\begin{prop}\label{Prop:Schanuel's Lemma resolution}
Suppose we have the following $\mathcal{M}$-injective resolutions of $E$, with the factorisation of each admissible morphism included:
\[\begin{tikzpicture}[auto]\matrix[column sep={0.5cm}, row sep={1.3cm,between origins}]{
	\node(E) {$E$}; &{} & \node(0)  {$I^0$};&{} & \node(dots)  {$\cdots$};&{} & \node(1)  {$I^{n-1}$};&{} & \node(m) {$I^n$};&{} & \node(F) {$\cdots$};\\
	\node(not1){ };&\node(E2) {$G^0$};&{}&\node(K1) {$G^1$};&{}&\node(K1m)  {$G^{n-1}$};&{}&\node(Km) {$G^n$}; &{}&\node(Km1) {$G^{n+1}$};\\	};	
	
	\draw[->] (E) to node {} (0); \draw[->] (0) to node {} (dots);
	\draw[->] (dots)to node {} (1); \draw[->] (1) to node {} (m);
	\draw[->] (m) to node {} (F);
		\draw[->>] (E) to node [swap]{{\footnotesize{$\cong$}}} (E2);
		\draw[->>] (0) to node { } (K1);\draw[->>] (dots) to node { } (K1m); \draw[->>] (1) to node { } (Km);
	\draw[>->] (E2) to node { } (0);\draw[>->] (K1) to node { } (dots); \draw[>->] (K1m) to node { } (1); \draw[>->] (Km) to node { } (m);
		\draw[->>]  (m) to node { } (Km1);
	\draw[>->] (Km1) to node { } (F);
\end{tikzpicture} 	\]
and
\[\begin{tikzpicture}[auto]\matrix[column sep={0.5cm}, row sep={1.3cm,between origins}]{
		\node(E) {$E$}; &{} & \node(0)  {$J^0$};&{} & \node(dots)  {$\cdots$};&{} & \node(1)  {$J^{n-1}$};&{} & \node(m) {$J^n$};&{} & \node(F) {$\cdots$};\\
		\node(not1){ };&\node(E2) {$H^0$};&{}&\node(K1) {$H^1$};&{}&\node(K1m)  {$H^{n-1}$};&{}&\node(Km) {$H^n$}; &{}&\node(Km1) {$H^{n+1}$};\\	};	
	
	\draw[->] (E) to node {} (0); \draw[->] (0) to node {} (dots);
	\draw[->] (dots)to node {} (1); \draw[->] (1) to node {} (m);
	\draw[->] (m) to node {} (F);
	\draw[->>] (E) to node [swap]{{\footnotesize{$\cong$}}} (E2);
	\draw[->>] (0) to node { } (K1);\draw[->>] (dots) to node { } (K1m); \draw[->>] (1) to node { } (Km);
	\draw[>->] (E2) to node { } (0);\draw[>->] (K1) to node { } (dots); \draw[>->] (K1m) to node { } (1); \draw[>->] (Km) to node { } (m);
		\draw[->>]  (m) to node { } (Km1);
	\draw[>->] (Km1) to node { } (F);
\end{tikzpicture} 	\]
Then, for each $n\geq{1}$, we have isomorphisms
\[ I^0{\oplus}J^1{\oplus}I^2{\oplus}\cdots{\oplus}J^{2n-1}{\oplus}G^{2n}
\,{}\cong{}\,
   J^0{\oplus}I^1{\oplus}J^2{\oplus}\cdots{\oplus}I^{2n-1}{\oplus}H^{2n}
\]
and
\[ I^0{\oplus}J^1{\oplus}I^2{\oplus}\cdots{\oplus}J^{2n-1}{\oplus}I^{2n}{\oplus}H^{2n+1}
\,{}\cong{}\,
J^0{\oplus}I^1{\oplus}J^2{\oplus}\cdots{\oplus}I^{2n-1}{\oplus}J^{2n}{\oplus}G^{2n+1}.
\]
\begin{proof}
We prove this by induction.
For $n=1$, first note that Corollary~\ref{Corollary: Schanuel's iso},
applied to the diagram
\[
\begin{tikzpicture}[auto]
	\matrix[column sep={1.1cm}, row sep={1.3cm,between origins}, matrix of math nodes](m)
	{ G^0 &   I^0 & G^1 \\
	  H^0 &   J^0  & H^1 \\
	};
	\draw[>->] (m-1-1) to node{ } (m-1-2); \draw[->>] (m-1-2) to node{ } (m-1-3);
	\draw[>->] (m-2-1) to node{ } (m-2-2); \draw[->>] (m-2-2) to node{ } (m-2-3);
	\draw[->]  (m-1-1) to node[swap]{{\footnotesize{$\cong$}}} (m-2-1);
\end{tikzpicture}
\]
gives
$I^0{\oplus}H^1 \cong J^0{\oplus}G^1$.
By Proposition~\ref{Prop: Direct sum kercokpair},
there is a diagram of the form
\[
\begin{tikzpicture}[auto]
	\matrix[column sep={1.1cm}, row sep={1.3cm,between origins}, matrix of math nodes](m)
	{ I^0{\oplus}H^1 &   I^0{\oplus}J^1\vphantom{H^2} & \vphantom{I^0{\oplus}J^1}H^2 \\
		J^0{\oplus}G^1 &   J^0{\oplus}I^1\vphantom{G^2}  & \vphantom{J^0{\oplus}I^1}G^2 \\
	};
	\draw[>->] (m-1-1) to node{ } (m-1-2); \draw[->>] (m-1-2) to node{ } (m-1-3);
	\draw[>->] (m-2-1) to node{ } (m-2-2); \draw[->>] (m-2-2) to node{ } (m-2-3);
	\draw[->]  (m-1-1) to node[swap]{{\footnotesize{$\cong$}}} (m-2-1);
\end{tikzpicture}
\]
and Corollary~\ref{Corollary: Schanuel's iso} gives $I^0{\oplus}J^1{\oplus}G^2  \cong J^0{\oplus}I^1{\oplus}{H^2}$.
To finish the proof for $n=1$,
we again apply Proposition~\ref{Prop: Direct sum kercokpair}
followed by Corollary~\ref{Corollary: Schanuel's iso}, to get a diagram
\[
\begin{tikzpicture}[auto]
	\matrix[column sep={1.1cm}, row sep={1.3cm,between origins}, matrix of math nodes](m)
	{ I^0{\oplus}J^1{\oplus}G^2\vphantom{G^3} &
	 I^0{\oplus}J^1{\oplus}I^2\vphantom{G^3} &
	 G^3 \\
	  J^0{\oplus}I^1{\oplus}{H^2}\vphantom{H^3} &
	  J^0{\oplus}I^1{\oplus}{J^2}\vphantom{H^3}  &
	  H^3 \\
	};
	\draw[>->] (m-1-1) to node{ } (m-1-2); \draw[->>] (m-1-2) to node{ } (m-1-3);
	\draw[>->] (m-2-1) to node{ } (m-2-2); \draw[->>] (m-2-2) to node{ } (m-2-3);
	\draw[->]  (m-1-1) to node[swap]{{\footnotesize{$\cong$}}} (m-2-1);
\end{tikzpicture}
\]
and an isomorphism
$I^0{\oplus}J^1{\oplus}I^2{\oplus}{H^3} \cong J^0{\oplus}I^1{\oplus}{J^2}{\oplus}{G^3}.$

\smallskip
Assume the result holds some $n\geq 1$.
By Proposition~\ref{Prop: Direct sum kercokpair},
there is a diagram of the form
\[
\begin{tikzpicture}[auto]
	\matrix[column sep={1.1cm}, row sep={1.3cm,between origins}, matrix of math nodes](m)
	{ I^0{\oplus}\cdots\oplus{I^{2n}}{\oplus}H^{2n+1} &
	  I^0{\oplus}\cdots\oplus{I^{2n}}{\oplus}J^{2n+1} &
	  H^{2(n+1)}
	  \\
	  J^0{\oplus}\cdots\oplus{J^{2n}}{\oplus}G^{2n+1} &
	  J^0{\oplus}\cdots\oplus{J^{2n}}{\oplus}I^{2n+1}\vphantom{G^{2(n+1)}} &
	  G^{2(n+1)} \\
	};
	\draw[>->] (m-1-1) to node{ } (m-1-2); \draw[->>] (m-1-2) to node{ } (m-1-3);
	\draw[>->] (m-2-1) to node{ } (m-2-2); \draw[->>] (m-2-2) to node{ } (m-2-3);
	\draw[->]  (m-1-1) to node[swap]{\footnotesize{$\cong$}} (m-2-1);
\end{tikzpicture}
\]
and Corollary~\ref{Corollary: Schanuel's iso} gives
\[ \begin{split}
	 &I^0{\oplus}J^1{\oplus}I^2{\oplus}\cdots{\oplus}J^{2(n+1)-1}{\oplus}G^{2(n+1)}\\
	 \cong \,&
	 J^0{\oplus}I^1{\oplus}J^2{\oplus}\cdots{\oplus}I^{2(n+1)-1}{\oplus}H^{2(n+1)}.
\end{split}
\]
One final application of Proposition~\ref{Prop: Direct sum kercokpair}
yields the following diagram:
\[
\begin{tikzpicture}[auto]
	\matrix[column sep={1.1cm}, row sep={1.5cm,between origins}, matrix of math nodes](m)
	{ {I^0{\oplus}J^1{\oplus}I^2{\oplus}\cdots{\oplus}J^{2n+1}{\oplus}G^{2(n+1)}}
	& J^0{\oplus}I^1{\oplus}J^2{\oplus}\cdots{\oplus}I^{2n+1}{\oplus}H^{2(n+1)}
		\\
		I^0{\oplus}J^1{\oplus}I^2{\oplus}\cdots{\oplus}J^{2n+1}{\oplus}I^{2(n+1)}
		& J^0{\oplus}I^1{\oplus}J^2{\oplus}\cdots{\oplus}I^{2n+1}{\oplus}J^{2(n+1)}
		\\
		G^{2(n+1)+1}
		& H^{2(n+1)+1}	\\
	};
	\draw[>->] (m-1-1) to node{ } (m-2-1); \draw[->>] (m-2-1) to node{ } (m-3-1);
	\draw[>->] (m-1-2) to node{ } (m-2-2); \draw[->>] (m-2-2) to node{ } (m-3-2);
	\draw[->]  (m-1-1) to node {{\footnotesize{$\cong$}}} (m-1-2);
\end{tikzpicture}
\]
By Corollary~\ref{Corollary: Schanuel's iso},
\[ \begin{split}
	& I^0{\oplus}J^1{\oplus}I^2{\oplus}\cdots{\oplus}J^{2(n+1)-1}{\oplus}I^{2(n+1)}{\oplus}H^{2(n+1)+1}\\
\cong\,&
J^0{\oplus}I^1{\oplus}J^2{\oplus}\cdots{\oplus}I^{2(n+1)-1}{\oplus}J^{2(n+1)}{\oplus}G^{2(n+1)+1}
\end{split}
\]
as required.
\end{proof}
\end{prop}

We can now prove the Injective Dimension Theorem.

\goodbreak
\begin{theorem}\label{Theorem: Injective dimension theorem}
Let $\mathcal{M}$ be the class of admissible monomorphisms in an exact category $(\mathcal{A},\Ex)$.
Suppose $\mathcal{A}$ has enough $\mathcal{M}$-injectives. The following are equivalent for $n\geq{1}$ and every $E\in\mathcal{A}$.
\begin{enumerate}[label=\upshape(\roman*)]
\item If there is an exact sequence of admissible morphisms
\begin{equation}\label{Diagram: injective dimension theorem diagram 1}
	\begin{tikzpicture}[auto, baseline=(current  bounding  box.center)]
				\matrix(m)[matrix of math nodes, column sep=3em]{
				{E}\vphantom{I^0} & {I^0} &	{\cdots}\vphantom{I^0} &   {I^{n-1}} & {F}\vphantom{I^0}
				\\};
				\draw[>->] (m-1-1) to node { } (m-1-2); \draw[->] (m-1-2) to node { } (m-1-3); \draw[->] (m-1-3) to node { } (m-1-4); \draw[->>] (m-1-4) to node { } (m-1-5);
	\end{tikzpicture}
		\end{equation}
		with each $I^{m}$, $0\leq m\leq n-1$ injective, then $F$ must be injective;
		\item There is an exact sequence of admissible morphisms	
		\begin{equation}\label{Diagram: injective dimension theorem diagram 2}\begin{tikzpicture}[auto, baseline=(current  bounding  box.center)]
				\matrix(m)[matrix of math nodes, column sep=3em]{
					{E}\vphantom{I^0} & {I^0} &	{\cdots}\vphantom{I^0} &   {I^{n-1}} & {I^n}\vphantom{I^0} \\};
				\draw[>->] (m-1-1) to node { } (m-1-2); \draw[->] (m-1-2) to node { } (m-1-3); \draw[->] (m-1-3) to node { } (m-1-4); \draw[->>] (m-1-4) to node { } (m-1-5);
			\end{tikzpicture}
		\end{equation}
		with each $I^{m}$, $0\leq m\leq n$ injective.
	\end{enumerate}	
\begin{proof}
Let $E\in\mathcal{A}$.
First we show (i) implies (ii). As $\mathcal{A}$ has enough $\mathcal{M}$-injectives, we can build an $\mathcal{M}$-injective resolution of $E$:
\[
\begin{tikzpicture}[auto]\matrix[column sep={0.5cm}, row sep={1.3cm,between origins}]{
		\node(E) {$E$}; &{} & \node(0)  {$J^0$};&{} & \node(dots)  {$\cdots$};&{} & \node(1)  {$J^{n-1}$};&{} & \node(m) {$J^n$};&{} & \node(F) {$\cdots$};\\
		\node(not1){ };&\node(E2) {$G^0$};&{}&\node(K1) {$G^1$};&{}&\node(K1m)  {$G^{n-1}$};&{}&\node(Km) {$G^n$}; &{}&{}\\	};	
	
	\draw[->] (E) to node {} (0); \draw[->] (0) to node {} (dots);
	\draw[->] (dots)to node {} (1); \draw[->] (1) to node {} (m);
	\draw[->] (m) to node {} (F);
	\draw[->>] (E) to node [swap]{{\footnotesize{$\cong$}}} (E2);
	\draw[->>] (0) to node { } (K1);\draw[->>] (dots) to node { } (K1m); \draw[->>] (1) to node { } (Km);
	\draw[>->] (E2) to node { } (0);\draw[>->] (K1) to node { } (dots); \draw[>->] (K1m) to node { } (1); \draw[>->] (Km) to node { } (m);
\end{tikzpicture} 	\]
Relabel $J^k$ as $I^k$ for all $0\leq k \leq n-1$ and $G^n$ as $I^n$,
this gives an exact sequence as in Diagram~\eqref{Diagram: injective dimension theorem diagram 2},
and $I^n$ must be $\mathcal{M}$-injective, by condition (i).

Now suppose that condition~(ii) holds.
There must exist an injective resolution of $E$ of the form
\[
\begin{tikzpicture}[auto]\matrix[column sep={0.5cm}, row sep={1.3cm,between origins}]{
		\node(E) {$E$}; &{} & \node(0)  {$J^0$};&{} & \node(dots)  {$\cdots$};&{} & \node(1)  {$J^{n-1}$};&{} & \node(m) {$J^n$};&{} & \node(F) {$\cdots$};\\
		\node(not1){ };&\node(E2) {$H^0$};&{}&\node(K1) {$H^1$};&{}&\node(K1m)  {$H^{n-1}$};&{}&\node(Km) {$J^n$}; &{}&{}\\	};	
	
	\draw[->] (E) to node {} (0); \draw[->] (0) to node {} (dots);
	\draw[->] (dots)to node {} (1); \draw[->] (1) to node {} (m);
	\draw[->] (m) to node {} (F);
	\draw[->>] (E) to node [swap]{{\footnotesize{$\cong$}}} (E2);
	\draw[->>] (0) to node { } (K1);\draw[->>] (dots) to node { } (K1m); \draw[->>] (1) to node { } (Km);
	\draw[>->] (E2) to node { } (0);\draw[>->] (K1) to node { } (dots); \draw[>->] (K1m) to node { } (1);
	\draw[>->] (Km) to node [swap]{{\footnotesize{$\cong$}}} (m);
\end{tikzpicture} 	\]
and for any exact sequence as in Diagram~\eqref{Diagram: injective dimension theorem diagram 1}, with each $I^n$ injective,
there exists an injective resolution
\[
\begin{tikzpicture}[auto]\matrix[column sep={0.5cm}, row sep={1.3cm,between origins}]{
		\node(E) {$E$}; &{} & \node(0)  {$I^0$};&{} & \node(dots)  {$\cdots$};&{} & \node(1)  {$I^{n-1}$};&{} & \node(m) {$I^n$};&{} & \node(F) {$\cdots$};\\
		\node(not1){ };&\node(E2) {$G^0$};&{}&\node(K1) {$G^1$};&{}&\node(K1m)  {$G^{n-1}$};&{}&\node(Km) {$G^n$}; &{}&{}\\	};	
	
	\draw[->] (E) to node {} (0); \draw[->] (0) to node {} (dots);
	\draw[->] (dots)to node {} (1); \draw[->] (1) to node {} (m);
	\draw[->] (m) to node {} (F);
	\draw[->>] (E) to node [swap]{{\footnotesize{$\cong$}}} (E2);
	\draw[->>] (0) to node { } (K1);\draw[->>] (dots) to node { } (K1m); \draw[->>] (1) to node { } (Km);
	\draw[>->] (E2) to node { } (0);\draw[>->] (K1) to node { } (dots); \draw[>->] (K1m) to node { } (1);
	\draw[>->] (Km) to node [swap]{ } (m);
\end{tikzpicture}
\]
with $G^n=F$.
By Proposition~\ref{Prop:Schanuel's Lemma resolution},
there exists a kernel-cokernel pair
\[
	\begin{tikzpicture}[auto, baseline=(current  bounding  box.center)]
		\node(F) {$F$};
		\node[right=11mm of F](I) {$I$};
		\node[right=11mm of I](G) {$G$};
		\draw[>->] (F) to node {\footnotesize{$\mu$}} (I);
		\draw[->>] (I) to node {\footnotesize{$\pi$}} (G);
	\end{tikzpicture}
\]
and a morphism $\widetilde{\mu}\in\Morunder{I}{F}{\mathcal{A}}$ such that
$\widetilde{\mu}\circ\mu=\idntyof{F}$,
and $I$ is a finite product of $\mathcal{M}$-injective objects.
Then by Proposition~\ref{Prop: Direct sum kercokpair}, $I$ is injective
and $\widetilde{\mu}$ is a left inverse for $\mu$, hence, by Proposition~\ref{Prop: M-injective},
$F$ is $\mathcal{M}$-injective.
\end{proof}
\end{theorem}

\begin{definition}\label{Definition: Injective Dimension}
Let $\mathcal{M}$ be the class of admissible monomorphisms in an exact category $(\mathcal{A},\Ex)$.
We say $E\in\mathcal{A}$ has \textit{finite $\mathcal{M}$-injective dimension} if there exists an exact sequence of admissible morphisms as in
Diagram~\eqref{Diagram: injective dimension theorem diagram 2} with all $I^m$ $\mathcal{M}$-injective. If $E$ is of finite
$\mathcal{M}$-injective dimension we write $\Midim{E}=0$ if $E$ is $\mathcal{M}$-injective and $\Midim{E}=n$ if $E$ is
not $\mathcal{M}$-injective and $n$ is the smallest natural number
such that there exists an exact sequence of admissible morphisms
as in Diagram~\eqref{Diagram: injective dimension theorem diagram 2} where every $I^m$ is $\mathcal{M}$-injective.
If $E$ is not of finite $\mathcal{M}$-injective dimension, we write $\Midim{E}=\infty$.
	
 The \textit{global dimension} of the exact category $(\mathcal{A},\Ex)$ is
 \[
\sup\setst{\Midim{E}}{E\in\mathcal{A}}\in\NN_0\cup\{\infty\}.
 \]
\end{definition}

\begin{remark}
The $\mathcal{M}$-injective dimension
of an object $E$ in an exact category $(\mathcal{A},\Ex)$
can be obtained by examining any of its $\mathcal{M}$-injective resolutions.
Indeed, suppose the following is an $\mathcal{M}$-injective resolution of $E$
(with the factorisation of each admissible morphism included):
\[\begin{tikzpicture}[auto]\matrix[column sep={0.5cm}, row sep={1.3cm,between origins}]{
		\node(E) {$E$}; &{} & \node(0)  {$J^0$};&{} & \node(dots)  {$\cdots$};&{} & \node(1)  {$J^{n-1}$};&{} & \node(m) {$J^n$};&{} & \node(F) {$\cdots$};\\
		\node(not1){ };&\node(E2) {$G^0$};&{}&\node(G1) {$G^1$};&{}&\node(G1m)  {$G^{n-1}$};&{}&\node(Gm) {$G^n$}; &{}&{}\\	};	
	
	\draw[->] (E) to node {} (0); \draw[->] (0) to node {} (dots);
	\draw[->] (dots)to node {} (1); \draw[->] (1) to node {} (m);
	\draw[->] (m) to node {} (F);
	\draw[->>] (E) to node [swap]{{\footnotesize{$\cong$}}} (E2);
	\draw[->>] (0) to node { } (G1);\draw[->>] (dots) to node { } (G1m); \draw[->>] (1) to node { } (Gm);
	\draw[>->] (E2) to node { } (0);\draw[>->] (G1) to node { } (dots); \draw[>->] (G1m) to node { } (1); \draw[>->] (Gm) to node { } (m);
\end{tikzpicture} 	\]
Then, by Theorem~\ref{Theorem: Injective dimension theorem},
 $\Midim{E}\leq n$ if and only if $G^n$ is $\mathcal{M}$-injective.
\end{remark}
The original version of Schanuel's lemma is formulated for projective resolutions, see, e.g., \cite[Lemma~5.1]{LamLectures}
or \cite[Theorem~3.41]{MR0409590}. An analogous version using the epimorphisms in the class $\mathcal P$ can be obtained in any
exact category with exact structure $(\mathcal M,\mathcal P)$.

\medskip
\bibliography{refs_modules}

\end{document}